\documentclass[12pt]{article}

\usepackage{amsfonts}
\usepackage{mathtools}
\usepackage{amsmath}
\usepackage{amssymb}
\usepackage{amsthm}
\usepackage{enumerate}

\def\parti{\partial^{\text{in}}}
\def\parto{\partial^{\text{out}}}

\def \Z {\mathbb Z}

\def \R {\mathbb R}

\def\P{\mathbb P}

\def\Ph1{P^{(h_1)}}
\def\Ph2{P^{(h_2)}}

\def\b0{{\bf 0}}

\newtheorem{thm}{Theorem}[section]
\newtheorem{prop}[thm]{Proposition}
\newtheorem{lem}[thm]{Lemma}
\newtheorem{cor}[thm]{Corollary}

\theoremstyle{plain}

\newtheorem{defn}[thm]{Definition}

\theoremstyle{definition}

\newtheorem{lemproof}[thm]{Proof of Lemma}

\begin{document}

\title{A lower bound for point-to-point connection probabilities in critical percolation}

\author{J. van den Berg\footnotemark[2] \footnote{During part of this project JvdB was affiliated
as visiting professor at NYU Abu Dhabi} 
$\,$ and H. Don\footnotemark[4] \\
\small \\
  \small \footnotemark[2] CWI and VU University Amsterdam ; email: J.van.den.Berg@cwi.nl \\
  \small \footnotemark[4] Radboud University Nijmegen; email: h.don@math.ru.nl
  }
\maketitle

\begin{abstract}
Consider critical site percolation on $\Z^d$ with $d \geq 2$.
We prove a lower bound of order $n^{- d^2}$
for point-to-point connection probabilities,
where $n$ is the distance between the points.

Most of the work in our proof concerns
 a `construction'
which finally reduces the problem to a topological one.
This is then solved by applying a topological fact,
Lemma \ref{topo-fact} below, which follows from Brouwer's fixed point theorem.

Our bound improves the lower bound with exponent $2 d (d-1)$, used  by Cerf in 2015 \cite{C15} to obtain an {\em upper} bound for
the so-called two-arm probabilities. 
Apart from being of interest in itself, our result gives a small improvement of the
bound on the two-arm exponent found by Cerf.

\end{abstract}

\begin{section}{Introduction and statement of the main results}

Consider site percolation on $\Z^d$ with $d \geq 2$. Let $p_c = p_c(d)$ denote the critical probability. We use the notation
$\Lambda(n)$ for the box $[-n,n]^d$. Our main result is the following, Theorem \ref{Mainthm} below. Although the result is,
essentially, only relevant for dimensions $3$ - $6$, we formulate it for all $d \geq 2$ because the proof works for all 
these values of $d$ `simultaneously'. 

\begin{thm}\label{Mainthm}
There is a constant $c>0$, which depends only on the dimension $d$, such that for
all $x \in \Z^d,x\neq 0$,

\begin{equation}\label{eq-Mainthm}
\P_{p_c} (0 \leftrightarrow x \text{ in } \Lambda(4 \|x\|)) \geq \frac{c}{\|x\|^{ d^2}}.
\end{equation}
Here $\|x\| := \max_{1\leq i \leq} |x_i|$, and $\{0 \leftrightarrow x \text{ in } \Lambda(4 \|x\|)\}$ is
the event that there is an open path from $0$ to $x$ which
lies in $\Lambda(4 \|x\|)$.

\end{thm}

Let us first explain why this result is, essentially, only relevant for dimensions $3$ - $6$:
For `nice' two-dimensional lattices (including the square lattice) a much better exponent ($2/3$ or even smaller) than
the exponent $4$ from \eqref{eq-Mainthm} is known. For these lattices it is believed (and known for site percolation
on the triangular lattice from \cite{LSW02}) that for every $C > 1$, 
$\P_{p_c} (0 \leftrightarrow x \text{ in } \Lambda(C \|x\|)) \approx {\|x\|}^{-5/24}$.

Further, it has been proved for dimensions $d \geq 11$ and is strongly believed for $d \geq 7$,
that the `unrestricted' event in \eqref{eq-Mainthm} (i.e. without the restriction ``in $\Lambda(4 \|x\|)$") 
has probability of order $\|x\|^{-(d-2)}$, see \cite{Har08} and \cite{HH17}. Moreover,
it seems to be believed that for such $d$ a suitably restricted form of this event
(with a suitable constant instead of the factor $4$ in the l.h.s. of \eqref{eq-Mainthm}) also has probability of
order $\|x\|^{-(d-2)}$.

The methods used for the above mentioned results for $d=2$ don't work for $d \in \{3, \dots,6\}$, and presumably this also holds for the methods used
for $d\geq 11$.
To our knowledge, for dimension $d \in \{3, \dots, 6\}$ no version of
Theorem \ref{Mainthm} with exponent $\leq d^2$ exists in the literature,
even if we drop the restriction ``in $\Lambda(4 \|x\|)"$.
In fact, we don't know if the probabilities for the events with and without the restriction are really of different order.
We stated our result for the case with restriction because it is stronger, and because of its
applicability for other purposes, see the Remarks below.

\medskip
Theorem \ref{Mainthm} has the following implication:
\begin{cor}\label{Maincor}
There is a constant $c' >0$, which depends only on the dimension $d$, such that for all $n \geq 1$,
and all $x, y \in \Lambda(n)$,

\begin{equation}\label{mainthm-bnd}
\P_{p_c} (x \leftrightarrow y \text{ in } \Lambda(9 n)) \geq \frac{c'}{n^{d^2}}.
\end{equation}
\end{cor}

\noindent
Indeed, this follows readily from Theorem \ref{Mainthm} by observing that, for all
$x, y \in \Lambda(n)$, we have that $\|y - x\| \leq 2n$ and that
$\P_{p_c} (x \leftrightarrow y \text{ in } \Lambda(9 n)) \geq \P_{p_c} (x \leftrightarrow y \text{ in } x + \Lambda(8 n))$,
which by translation invariance is equal to $\P_{p_c} (0 \leftrightarrow y-x \text{ in }  \Lambda(8 n))$.

\medskip
Corollary \ref{Maincor} clearly improves (apart from the factor $9$ in the l.h.s.\ of \eqref{mainthm-bnd}, see Remark (ii) below) the following
result by Cerf, which is one of the motivations for our work:

\begin{lem}\label{cerf-lem6.1} {\em (Cerf, Lemma 6.1 in \cite{C15})} \\
There is a constant $\tilde c >0$, which depends only on the dimension $d$, such that for all $n \geq 1$,
and all $x, y \in \Lambda(n)$,

\begin{equation}\label{cerf-bnd}
\P_{p_c} (x \leftrightarrow y \text{ in } \Lambda(2 n)) \geq \frac{\tilde c}{n^{2 d (d-1)}}.
\end{equation}
\end{lem}

\medskip\noindent
{\bf Remarks:} 

(i) Lemma \ref{cerf-lem6.1} above was used in a clever way by Cerf in \cite{C15} to obtain new {\it upper} 
bounds for two-arm probabililities.
Cerf's work on two-arm probabililities is, in turn, mainly motivated by the famous 
conjecture that, for any dimension $\geq 2$, there is no infinite open cluster at $p_c$.
This conjecture, one of the main open problems in
percolation theory, has only been proved for dimension two and for high
dimensions. 

(ii) The precise factor $2$ in the expression $\Lambda(2 n)$ in \eqref{cerf-bnd} is not essential for Cerf's applications
mentioned in Remark (i):
a bigger constant, which may even depend on the dimension, would work as well (with tiny, straightforward modifications
of Cerf's arguments). 
This is why we have not seriously tried to reduce the factor $9$ in the expression $\Lambda(9 n)$ in Corollary \ref{Maincor}.

(iii) What {\it does} matter for Cerf's applications is an improvement of the 
power $2 d(d-1)$ in the r.h.s.\ of \eqref{cerf-bnd}.
By a simple adaptation of a step in Cerf's proof of Lemma \ref{cerf-lem6.1}, the
exponent $2 d (d-1)$  can be replaced by $(2d - 1)(d-1)$. Our Corollary \ref{Maincor} states that it can be
replaced by $d^2$.
So, for instance, for $d=3$ the exponents given by Cerf, by the small adaptation of Cerf's argument mentioned above,
and by our Corollary \ref{Maincor}, are $12$, $10$ and $9$ respectively; and for $d=4$ they are $24$, $21$ and $16$ respectively.
Following step by step the computations in Sections 7-9 in Cerf's paper, shows that our improvement of Lemma \ref{cerf-lem6.1} 
also provides a small improvement of Cerf's lower bound for the two-arm exponent.

\begin{subsection}{Main ideas in the proof of Theorem \ref{Mainthm}}

We first present, somewhat informally, Cerf's proof of Lemma \ref{cerf-lem6.1} (and the proof of the small
adaptation mentioned in Remark (iii) above).

To focus on the main idea we ignore here the restriction ``in $\Lambda(2n)$" in Lemma \ref{cerf-lem6.1}.
The key ingredient is a result, Corollary \ref{cor-hamm} below, which goes back to work by Hammersley.
Applied to the special case where $\Gamma$ is a box centered at $0$, and using symmetry, this result gives that there is a $C = C(d)$ such that
for every $n$ there is a `special' site $v^n$ on the boundary of $\Lambda(n)$ with $v^n_1 = n$ and such that 
$$\P_{p_c}(0 \leftrightarrow v^n) \geq \frac{C}{n^{d-1}}.$$

Now let $x$ be a vertex. 
By symmetry we may assume that its coordinates $x_i$, $i=1, \dots ,d$, are non-negative. 
For simplicity we also
assume that they are even.
Let $x(i)$ be the vertex of which the $i$th coordinate is equal to $x_i$ and all other coordinates are $0$.

We have by translation invariance and FKG that $\P_{p_c} (0 \leftrightarrow x)$ is larger than or equal to
$$\prod_{i=1}^d \P_{p_c}(0 \leftrightarrow x(i)).$$

Using symmetry, the $i$th factor in this product is equal to 
$\P_{p_c}(0 \leftrightarrow (x_i, 0, 0, \cdots,0))$,
and hence, by FKG (and again by symmetry), at least 
$$(\P_{p_c}(0 \leftrightarrow v^{x_i/2}))^2,$$
which as we saw is $\geq C^2\cdot{\|x\|}^{-2(d-1)}$. 
Hence the product of the $d$ factors is at least of order
$\|x\|^{-2 d (d-1)}$. This gives essentially (apart from some relatively straightforward issues) Lemma \ref{cerf-lem6.1}.

The small adaptation mentioned before comes from an observation concerning the first step of the argument above:
Note that $\P_{p_c}(0 \leftrightarrow x)$ is also larger than or equal to 
$$\P_{p_c}(0 \leftrightarrow v^{x_1}) \, \P_{p_c}(v^{x_1} \leftrightarrow x).$$
The first factor is at least $C\cdot \|x\|^{-(d-1)}$. And the second factor can be written as
$\P_{p_c}(0 \leftrightarrow v^{x_1} - x)$ to which we can apply Cerf's argument mentioned above. However, since 
$v^{x_1} - x$ has first coordinate $0$, and hence has at most $d-1$ non-zero coordinates, that argument gives a 
lower bound of order $\|x\|^{-2(d-1)(d-1)}$.
Hence, the product of the two factors is at least of order $\|x\|^{-(2 d -1)(d-1)}$ in agreement with Remark (iii).

 So, roughly speaking, the explanation of Lemma \ref{cerf-lem6.1} is that $x$ can be written as the sum of $2 d$ 
(and in fact, as the above adaptation showed, $2 d - 1$) special points, where each special
point `costs' a factor of order $\|x\|^{-(d-1)}$.

Our proof of Theorem \ref{Mainthm} (and thus of Corollary \ref{Maincor}) also uses the idea of certain special points, 
which will be called `good' points.
Our notion of good points is weaker than that of the special points mentioned above, in the sense that
they `cost' a factor of order $\|x\|^{-d}$. However, we will prove that, roughly speaking, each point is the sum of `just' $d$ good points, which
together with the previous statement gives the exponent $d^2$ in Theorem \ref{Mainthm}.
To prove this we first show, again using Corollary \ref{cor-hamm} (but now with general $\Gamma$, not only boxes)
the existence of suitable paths of good points and then turn the problem into a topological issue.
This is then finally solved by applying the `topological fact' Lemma \ref{topo-fact}.

\end{subsection}

\end{section}

\begin{section}{Proof of Theorem \ref{Mainthm}}
\begin{subsection}{Reformulation of Theorem \ref{Mainthm}}
Before we start, we reformulate Theorem \ref{Mainthm} as follows:

\begin{prop}\label{Mainprop}
There is a constant $c>0$, which depends only on the dimension $d$, such that for all $n \geq 1$,
and all $x \in \Lambda(n)$,

\begin{equation}\label{eq-Mainprop}
\P_{p_c} (0 \leftrightarrow x \text{ in } \Lambda(4 n)) \geq \frac{c}{n^{ d^2}}.
\end{equation}
\end{prop}

This proposition is trivially equivalent to Theorem \ref{Mainthm}. This reformulation is less `compact' than that of 
Theorem \ref{Mainthm} but has the advantage that it is more natural with respect to the approach in our proof (where, as we will
see, we first fix an $n$ and then distinguish between `good' and other points in $\Lambda(n)$.

\end{subsection}
\begin{subsection}{Preliminaries}
First we introduce some notation and definitions: 
For two vertices $v, w$, the notation $v\sim w$ means that $v$ and $w$ are neighbours, i.e. that
$\|v - w\|_1 := \sum_{i=1}^d |v_i - w_i| = 1$. \\

If $V$ and $W$ are sets of vertices, a {\it path} from $V$ to $W$ is a sequence of vertices
$v(1), \dots, v(k)$ with $v(1) \in V$, $v(k) \in W$ and $\|v(i+1) - v(i)\|_1 \leq 1$ for all
$i = 1, \dots,k-1$. (Note that we allow consecutive vertices to be equal; this is done for convenience later in this paper.)

\smallskip
If $W$ is a set of vertices, we define
\begin{align}\nonumber
 & \parto \, W  = \{v \in W^c \, : \, \exists w \in W \text{ s.t. } w \sim v\}, \,\,\, \text{ and } \\ \nonumber
 & \parti \, W = \{w \in W \, : \, \exists v \in W^c \text{ s.t. } w \sim v\}.
\end{align}

\noindent 

\medskip
We will use a result, Lemma \ref{hammersley} below, which goes back to work by Hammersley in the late fifties. 
Hammersley only proved the special case (used by Cerf) where $\Gamma$ is of the form $\Lambda(k)$.
The proof of the general case (which we use) is very
similar, it is given in Section 2 of \cite{DC-T16}.

\begin{lem}\label{hammersley} (Hammersley \cite{H57}, Duminil-Copin and Tassion \cite{DC-T16}.)
Let $\Gamma$ be a finite, connected set of vertices containing $0$. Then

\begin{equation}\label{eq-hammersley}
\sum_{x \in \parto \, \Gamma} \P_{p_c}(0 \leftrightarrow x \text{ `in' } \, \Gamma) \geq 1,
\end{equation}
where $\{0 \leftrightarrow x \text{ `in' } \Gamma\}$ is the event that there is an open path from $0$ to $x$ of which all vertices,
except $x$ itself, are in $\Gamma$.

\end{lem}

Clearly, this lemma has the following consequence.

\begin{cor}\label{cor-hamm}
For every finite, connected set of vertices $\Gamma$ containing $0$,
there is an $x \in \parto \, \Gamma$ such that

\begin{equation}\label{eq-hamm}
\P_{p_c}(0 \leftrightarrow x \text{ `in' } \Gamma) \geq \frac{1}{|\parto \, \Gamma|}.
\end{equation}
\end{cor}

\end{subsection}

\begin{subsection}{The set of `good' vertices and its properties}

Since, throughout the proof, the percolation parameter $p$ is equal to $p_c$, we simply write
$\P$ instead of $\P_{p_c}$.

\smallskip
Let $n\geq 1$ be fixed, and consider the box $\Lambda(n)$.

\begin{defn}
A vertex $v \in \Lambda(n)$ is called {\em good} if $v = 0$ or
$$\P(0 \leftrightarrow v \text{ in } \Lambda(n)) \geq \frac{1}{|\Lambda(n)|}.$$
The set of all good vertices is denoted by $G$.
\end{defn}

\begin{lem}\label{G-props} The set $G$ has the following properties:
\begin{enumerate}[(a)]
\item $G$ is invariant under the automorphisms of $\Lambda(n)$. 
\item For each $i =1, \dots, d$, there is a path $\pi^{(i)} \subset G \cap [0,n]^d$ which starts in $0$ and ends
in $[0,n]^{i-1} \times \{n\} \times [0,n]^{d-i}$.
\end{enumerate}
\end{lem}

\addtocounter{thm}{-1}
\begin{lemproof}\ \\
\noindent Part (a) is obvious. \\
As to part (b), we define, inductively, a finite sequence $\Gamma_1 \subset \Gamma_2 \subset \cdots$ of subsets of $\Lambda(n)$
as follows: $\Gamma_1 := \{0\}$. If $j \geq 1$ and $\Gamma_j$ is given, then, if
$\Gamma_j \cap \partial^{\text{in}}\Lambda(n)$ is nonempty we stop the procedure.
If it is empty, we can find, by Corollary \ref{cor-hamm},
a vertex $x \in \parto \, \Gamma_j \cap \Lambda(n)$ such that 
$$\P(0 \leftrightarrow x \text{ in } \Lambda(n)) \geq \P(0 \leftrightarrow x \text{ `in' } \Gamma_j) \geq
\frac{1}{|\parto \, \Gamma_j|} \geq \frac{1}{|\Lambda(n)|}.$$ 
In that case we take $\Gamma_{j+1} := \Gamma_j \cup \{x\}$.

Let $\Gamma$ denote the final set obtained at the end of the procedure (i.e. $\Gamma$ is the union of the $\Gamma_j$'s obtained by the procedure).
It is immediate from the procedure that all the vertices of $\Gamma$ are good, and that there is a path
$\pi = (x(1) = 0, \, x(2), \dots, x(m))$ in $\Gamma$ from
$0$ to $\partial^{\text{in}} \Lambda(n)$. Clearly, (where, for $y = (y_1, \dots, y_d) \in \R^d$, we 
use the notation $|y|$ for $(|y_1|, \dots, |y_d|)$),
$(|x(1)|, |x(2)|, \dots, |x(m)|)$ is a path in $[0,n]^d$, which starts in $0$ and of which the endpoint is 
in $[0,n]^{i-1} \times \{n\} \times [0,n]^{d-i}$ for some $i \in \{1, \dots,d\}$.
Moreover, by part (a) of the lemma, each vertex of this
path belongs to $G$. Hence, the property in part (b) of the lemma holds for some $i \in \{1, \dots,d\}$.
Using again part (a), the property follows for
{\it each} $i \in \{1, \dots,d\}$. \qed
\end{lemproof}

\begin{lem}\label{G-props-x}
For each $i =1, \dots, d$ there is a path $\tilde \pi^{(i)}$ of good points in the box $[0,n]^{i-1} \times [-n,n] \times [0,n]^{d-i}$, starting at a point
of which the $i$th coordinate is $-n$ and ending at a point of which the $i$th coordinate is $n$. (So $\tilde \pi^{(i)}$ traverses the box `in the long direction').

\end{lem}

\addtocounter{thm}{-1}
\begin{lemproof}\ \\
Fix an $i \in \{1, \dots,d\}$. For each vertex $v = (v_1, \dots, v_d)$, let the vertex $\tilde v = (\tilde v_1, \dots, \tilde v_d)$ be defined by:

$$ 
\tilde v_j = 
\begin{cases}
v_j, & j \neq i \\
- v_j, & j=i
\end{cases}
$$

\noindent
Note that, by Lemma \ref{G-props}(a), $v \in G$ iff $\tilde v \in G$.\\
Let $\pi^{(i)}$ be the path in the statement of Lemma \ref{G-props}(b),
and
write this path as $(v(1), v(2), \dots, v(k))$. Then, clearly, \\
$(\tilde v(k), \tilde v(k-1), \dots, \tilde v(2), \tilde v(1), v(1), v(2), \dots, v(k))$
is a path with the desired property. \qed
\end{lemproof}
\medskip
Now, for each $i = 1, 2, \dots, d$, we associate the path $\tilde \pi^{(i)}$ of Lemma \ref{G-props-x} with a continuous curve $c^{(i)}$,
simply by following the path along the edges of the graph at constant speed, chosen such that the total travel time along the path is $1$.
By Lemma \ref{G-props-x} we immediately have the following:

\begin{lem}\label{c-props}
For each $i =1, \dots, d$, the function 
$c^{(i)} \, : \, [0,1] \rightarrow \R^d$ is continuous and has the following properties:

\medskip\noindent
 (i) $\forall s \in [0,1] \,\, \exists t \in [0,1] \text{ s.t. } c^{(i)}(t) \in G \text{ and } \|c^{(i)}(t) - c^{(i)}(s)\|_1 \leq 1/2$.\\
(ii) For all $s \in [0,1]$, $\, -n = (c^{(i)}(0))_i \leq (c^{(i)}(s))_i \leq (c^{(i)}(1))_i = n$. \\
 (iii) $(c^{(i)}(s))_j \in [0,n]$ for all $s \in [0,1]$ and $j \neq i$.

\end{lem}

\end{subsection}

\subsection{Topological approach}

We use the curves $c^{(i)}, i = 1, \ldots, d$, introduced in the previous subsection, to define a
(as it turns out, useful) function $g \, : \, [0,1]^d \rightarrow \R^d$. First we define for each $t = (t_1,\ldots,t_d)\in [0,1]^d$ the matrix $C(t)\in\mathbb{R}^{d\times d}$ by
\begin{equation}
C(t) = \left[c^{(1)}(t_1) \quad c^{(2)}(t_2)\qquad \cdots\qquad c^{(d)}(t_d)\right].
\end{equation}
Next, we introduce functions $h_d:\mathbb{R}^d\rightarrow \mathbb{R}$ defined recursively for $d\geq 1$ by
\begin{equation}
\begin{array}{ll}
h_1(x) = |x| & \text{for}\ x\in \mathbb{R},\\
h_d(x) = |h_{d-1}(x_1,\ldots,x_{d-1})-x_d|  & \text{for}\ x \in \mathbb{R}^d,d\geq 2.
\end{array}
\end{equation}
By induction it easily follows that for all $d\geq 1$ and all $x\in\mathbb{R}^d$ there exist weights $a_j\in\left\{-1,1\right\}$ for $j=1,\ldots,d$ such that
\begin{equation}\label{h-prop}
0\leq h_d(x) = \sum_{j=1}^d a_j x_j\leq \Vert x\Vert \qquad\text{and}\qquad \forall k\leq d: \left|\sum_{j=1}^k a_jx_j\right| \leq \Vert x\Vert.
\end{equation}
Now we define $g:[0,1]^d\rightarrow\mathbb{R}^d$ by describing the $d$ coordinates of $g(t)$. 
\begin{defn}\label{g-def} For $t\in [0,1]^d$, let $C=C(t)$. For all $1\leq i\leq d$, we define 
\[
(g(t))_i =  C_{i,i} + h_{d-1}(C_{i,1},\ldots,C_{i,i-1},C_{i,i+1},\ldots,C_{i,d}).
\]
\end{defn}

\begin{lem}\label{g-props}
The function $g:[0,1]^d\rightarrow \mathbb{R}^d$ has the following properties:
\begin{enumerate}[(a)]
	\item\label{cont} $g$ is continuous;
	\item\label{rand} $(g(t))_i\leq 0$ if $t_i=0$\qquad and\qquad  $(g(t))_i\geq n$ if $t_i=1$;
	\item For each $x\in g([0,1]^d)$, there exist $z^{(1)},\ldots,z^{(d)}\in G$ such that
	\begin{enumerate}[(i)]
		\item $\sum_{j=1}^k z^{(j)} \in \Lambda(2n+\lceil\frac d 2\rceil)$ for $k=1,\ldots,d$;
		\item $\left\Vert x-\sum_{j=1}^d z^{(j)}\right\Vert_1\leq \frac d 2$.
	\end{enumerate} 
\end{enumerate}	
\end{lem}

\addtocounter{thm}{-1}
\begin{lemproof}\ \\
Part (\ref{cont}) follows because $g$ is a composition of continuous functions. Note that $C_{i,j}(t) = (c^{(j)}(t_j))_i$ so that $|C_{i,j}(t)|\leq n$ for all $i,j$ by Lemma \ref{c-props}(ii) and (iii). By (\ref{h-prop}), this implies that  
\[
0\leq h_{d-1}(C_{i,1},\ldots,C_{i,i-1},C_{i,i+1},\ldots,C_{i,d}) \leq n.
\]
Together with Lemma \ref{c-props}(ii), this proves part (\ref{rand}).

\noindent By definition of $g$ and by (\ref{h-prop}), for each $t\in [0,1]^d$ there exist weights $a_{ij}(t)\in\left\{-1,1\right\}$ such that
\begin{equation}\label{eq-gsum}
g(t) = \tilde C(t)\mathbf{1} = \sum_{j=1}^{d} \tilde c^{(j)}(t_j),
\end{equation}
where $\mathbf{1}$ denotes $(1, \dots,1)$ and $\tilde C(t)$ is the matrix with elements $\tilde C_{i,j}(t) = a_{ij}(t)C_{i,j}(t)$ and columns $\tilde c^{(j)}(t_j)$. For all $i$, the partial sums of (\ref{eq-gsum}) satisfy
\[
\biggl(\sum_{j=1}^k \tilde c^{(j)}(t_j)\biggr)_i = \left\{\begin{array}{ll}
\sum_{j=1}^k a_{ij}C_{i,j} & \text{if}\ k<i,\\
C_{ii}+\sum_{j=1,j\neq i}^k a_{ij}C_{i,j}\qquad &\text{if}\ k\geq i.
\end{array}\right.
\]
Applying (\ref{h-prop}) gives that all partial sums are in $\Lambda(2n)$ for all $t$. 

Since Lemma \ref{c-props}(i) applies to the columns of $C(t)$, it also applies to the columns of $\tilde C(t)$ by Lemma \ref{G-props}(a). Consequently, there exist $s = (s_1,\ldots,s_d) \in[0,1]^d$ such that $z^{(j)}  := \tilde c^{(j)}(s_j)\in G$ for all $j$ and such that for all $k\leq d$
\[
\biggl\Vert \sum_{j=1}^{k} \tilde c^{(j)}(t_j)-\sum_{j=1}^k z^{(j)}\biggr\Vert_1 \leq \sum_{j=1}^k \left\Vert\tilde c^{(j)}(t_j)-z^{(j)}\right\Vert_1 \leq \frac k 2\leq \frac d 2,
\]
completing the proof of part (c).\qed
\end{lemproof}

\begin{lem}\label{MainProp-for-g}
There exists a constant $\tilde c > 0$ (which depends on $d$ only) such that the inequality in Proposition \ref{Mainprop} 
holds for all lattice points $ x \in g([0,1]^d) \cap \Lambda(n)$. 

\end{lem}

\addtocounter{thm}{-1}
\begin{lemproof}\ \\
Let $x$ be a lattice point in $g([0,1]^d) \cap \Lambda(n)$. If $n\leq d/2$, then 
\[
\P(0 \leftrightarrow x \text{ in } \Lambda(4n)) \geq p_c^{\|x\|_1} \geq p_c^{d^2/2}.
\]
So we can assume that $3n+d/2 < 4n$. Choose $z^{(1)},\ldots,z^{(d)}$ as in Lemma \ref{g-props}(c). We get that the event 
$\{0 \leftrightarrow x \text{ in } \Lambda(4n)\}$ occurs if each of  the events
$$\left\{0 \stackrel{ \text{ in } \Lambda(n)} {\longleftrightarrow} z^{(1)}\right\},$$ 
$$\bigcap_{k=1}^{d-1} \biggl\{\sum_{i=1}^k z^{(i)} 
\longleftrightarrow \sum_{i=1}^{k+1} z^{(i)} \text{ in } \sum_{i=1}^k z^{(i)} + \Lambda(n) \biggr\}, \text{ and }$$ 
$$ \biggl\{\sum_{i=1}^d z^{(i)}
\leftrightarrow x \text{ in }\sum_{i=1}^d z^{(i)} + \Lambda(d)\biggr\}$$
occurs.
Hence, by FKG and using that each $z^{(i)}$ is good, 
$$ \P(0 \leftrightarrow x \text{ in } \Lambda(4n)) \geq \left(\frac{1}{(2n+1)^d}\right)^d \, p_c^{d/2}\geq \frac{p_c^{d/2}}{3^{d^2}}\cdot \frac 1 {n^{d^2}}.$$
This completes the proof of Lemma \ref{MainProp-for-g}. \qed 
\end{lemproof}

\begin{subsection}{Completion of the proof of Theorem \ref{Mainthm}} 
By Lemma \ref{MainProp-for-g} it is clear that the result of Proposition \ref{Mainprop}, and hence Theorem \ref{Mainthm}, 
follows if 
\begin{equation}\label{eq-semi-final}
g([0,1]^d) \supset [0,n]^d.
\end{equation}
Indeed, if $g([0,1]^d) \supset [0,n]^d$, 
then, by Lemma \ref{MainProp-for-g}, 
the bound \eqref{eq-Mainthm} holds for all vertices in $[0,n]^d$, and hence, by symmetry (part (a) of Lemma \ref{G-props}),
also for all vertices in $\Lambda(n)$.

Finally, \eqref{eq-semi-final} follows by applying the following `topological fact', Lemma \ref{topo-fact}, 
to the function $\frac{g(\cdot)}{n}$, and by noting that Lemma \ref{g-props} (a) and (b) guarantees that this function
satisfies the conditions of Lemma \ref{topo-fact}.
We don't know if this lemma is, more or less explicitly, in the literature. We are grateful to Lex Schrijver for
providing his proof.

\begin{lem}\label{topo-fact}
Let $f \, : \, [0,1]^d \rightarrow \R^d$ be a continuous function with the following properties.

\begin{equation}\label{eq-top-0}
\text{For all } i \in \{1, \dots, d\} \text{ and all } x \in [0,1]^d \text{ with } x_i = 0, \,\, (f(x))_i \leq 0,
\end{equation}

\begin{equation}\label{eq-top-1}
\text{For all } i \in \{1, \dots, d\} \text{ and all } x \in [0,1]^d \text{ with } x_i = 1, \,\, (f(x))_i \geq 1,
\end{equation}

Then $f([0,1]^d) \supset [0,1]^d$.
\end{lem}

\addtocounter{thm}{-1}
\begin{lemproof}\ [Lex Schrijver, private communication]\\
Since $f([0,1]^d)$ is compact, it is sufficient to prove that $f([0,1]^d) \supset (0,1)^d$.
Further, w.l.o.g. we may assume that $f([0,1]^d) \subset [0,1]^d$ and that
the inequalities at the end of \eqref{eq-top-0} and \eqref{eq-top-1} are {\it equalities}.
To see this, replace $f$ by the function $\tilde f$ defined by

$$ 
(\tilde f(x))_i =
\begin{cases}
1, & (f(x))_i \geq 1 \\
0, & (f(x))_i \leq 0 \\
(f(x))_i, & \text{otherwise}
\end{cases}
$$
and note that $\tilde f([0,1]^d) \cap (0,1)^d = f([0,1]^d) \cap (0,1)^d$.

\smallskip
So, from now on, we assume that $f: [0,1]^d \rightarrow [0,1]^d$ is continuous and has the property that: \\
For all  $i \in \{1, \dots, d\}$ and all $x \in [0,1]^d$  with  $x_i = 0$, $(f(x))_i = 0$, and \\
for all  $i \in \{1, \dots, d\}$ and all $x \in [0,1]^d$  with $x_i = 1$,  $(f(x))_i = 1$. \\
Note that this implies in particular that, for all $x \in \partial [0,1]^d$, $f(x) \in \partial [0,1]^d$ and
\begin{equation}\label{eq-forb}
f(x) \neq {\bf 1} - x.
\end{equation}

Now suppose there is an $a \in (0,1)^d)$ which is not in $f([0,1]^d)$. We will see that this leads to a contradiction.
Fix such an $a$ and let $\pi$ be the projection from $a$ on $\partial[0,1]^d$. (More precisely, for every $x \neq a$,
$\pi(x)$ is the unique intersection of the half-line $\{a + \lambda (x-a) \, : \, \lambda >0\}$ with $\partial[0,1]^d$). 
Finally, let $\psi(x) := {\bf 1} - x$.

Now consider the function $\psi \circ \pi \circ f \, : \, [0,1]^d \rightarrow \partial [0,1]^d$.
Note that it is well-defined (by the special property of $a$) and continuous. Hence, by Brouwer's fixed-point theorem, it has a fixed point,
which we denote by $b$. It is clear that $b \in \partial [0,1]^d$. So we have 

$$ b = (\psi \circ \pi \circ f)(b) = \psi(f(b)),$$
and hence $\psi(b) = f(b)$. However, this last equality violates \eqref{eq-forb}, so we have indeed obtained a contradiction.
\qed
\end{lemproof}

\end{subsection}

\end{section}

\bigskip\noindent
{\large\bf Acknowledgments}\\
We are grateful to Lex Schrijver for his proof of Lemma \ref{topo-fact}. We also thank Steffen Sagave and Ben Moonen for their
help with a topological question which came up during an earlier, somewhat different attempt to prove Theorem \ref{Mainthm}.
We thank Wouter Cames van Batenburg for very stimulating discussions, and Remco van der Hofstad for interesting
information concerning two-point probabilities in high dimensions.

\end{document}